\documentclass[letterpaper,11pt]{article}

\usepackage{amsmath,amsthm}
\usepackage{amsfonts}
\usepackage{latexsym}
\usepackage{amssymb,bm}
\usepackage{mathrsfs}
\usepackage{mathtools}
\usepackage[page]{appendix}

\newtheorem{theorem}{Theorem}[section]
\newtheorem{lemma}[theorem]{Lemma}

\newtheorem{proposition}[theorem]{Proposition}
\DeclareMathOperator*{\argmin}{argmin}

\newcommand{\inner}[2]{\langle #1,#2\rangle}
\newcommand{\Norm}[1]{\left\|{#1}\right\|}
\newcommand{\norm}[1]{\|{#1}\|}
\newcommand\set[1]{\{#1\}}

\newcommand{\R}{\mathbb{R}}
\newcommand{\tos}{\rightrightarrows}

\newcommand{\N}{\mathbb{N}}

\title{Complexity of the relaxed Hybrid Proximal-Extragradient method under the large-step condition}

 \author{ Benar F. Svaiter\thanks{ IMPA, Estrada Dona
     Castorina 110, 22460-320 Rio de Janeiro, Brazil ({\tt benar@impa.br}).}
   \hspace{.5em}\thanks{This work was
   Partially supported by CNPq grants
   302962/2011-5, %
   474996/2013-1, %
   FAPERJ grant
   201.584/2014 (Cientista d0 Nosso Estado)
   and by PRONEX-Optimization.}}

\begin{document}

\maketitle

\begin{abstract}
In this note we review the iteration-complexity of a relaxed Hybrid-Proximal
Extragradient Method under the large step condition. We also derive some
useful proprieties of this method.
\\
\\ 
  2000 Mathematics Subject Classification: 90C60, 90C25, 47H05.
 \\
 \\
 Key words:
 monotone inclusion,
 hybrid proximal extragradient method, large step condition, complexity, ergodic convergence, relaxation, maximal monotone operator.
\end{abstract}

\pagestyle{plain}

\subsection*{Introduction}

In this note we review Rockafellar's Proximal Point Method and the
Relaxed Hybrid Proximal-Extragradient (r-HPE) Method. We also some
useful properties of the (r-HPE) and analyze its complexity under the
large-step condition.
All the presented results pertaining the r-HPE were essentially proved
in~\cite{pre-print-impa}. The unique exception are the first inequalities
in  item~\ref{it:5-a} of
Lemma~\ref{lm:a} and in item~\ref{it:4-b} of Proposition~\ref{pr:a1}.

\section{Maximal monotone operators, the monotone inclusion problem,
and Rockafellar's Proximal Point Method}

Let $ H $ be a Hilbert space with inner product $\inner{\cdot}{\cdot}$ and
associated norm $\norm\cdot$. A \emph{point-to-set
  operator} in $H$, $T: H \tos H $, is a relation $T\subset  H \times H $ and
\[
T(z)=\{v\;|\; (z,v)\in T\},\qquad z\in  H .
\]
The \emph{inverse} of $T$ is $T^{-1}: H \tos H $,
$T^{-1}=\{(v,z)\;|\; (z,v)\in T\}$. The \emph{domain} and the range of
$T$ are, respectively,
\[
D(T)=\set{z\;|\; T(z)\neq\emptyset},\quad R(T)=\set{v\;|\;\exists z\in
   H ,\; v\in T(z)}.
\]
When $T(z)$ is a singleton for all $z\in D(T)$ it is usual to identify $T$
with the map $D(T)\ni z \mapsto v\in H $ where $T(z)=\set{v}$.
If $T_1,T_2: H \tos H $ and $\lambda\in\R$, then $T_1+T_2: H \tos H $ and
$\lambda T_1: H \tos H $ are defined as
\[
(T_1+T_2)(z)=\set{v_1+v_2\;|\; v_1\in T_1(z),\;v_2\in T_2(z)},\qquad
(\lambda T_1)(z)=\set{\lambda v\;|\; v\in T_1(z)}.
\]

A point-to-set operator $T: H \tos H $ is \emph{monotone} if
\[
\inner{z-z'}{v-v'}\geq 0\qquad \forall (z,v),(z',v')\in T
\]
and it is \emph{maximal monotone} if it is a maximal element on the family of
monotone point-to-set operators in $ H $ with respect to the partial order of
set inclusion.
Minty's theorem~\cite{MR0169064} states that if $T$ is maximal monotone
and $\lambda>0$, then the \emph{proximal map} $(\lambda T+I)^{-1}$ is a
point-to-point non-expansive operator with domain $ H $.

The \emph{monotone inclusion problem} is: given $T: H \tos H $ maximal
monotone, find $z$ such that
\begin{align}
  \label{eq:g.mip}
  0\in T(z).  
  \end{align}
  Rockafellar's Proximal Point Method~\cite{Rock:ppa} (hereafter PPM)
generates, for any starting $z_0\in H $, a sequence $(z_k)$ by the
approximate rule
\begin{align*}
  z_k\approx (\lambda_k T+I)^{-1}z_{k-1},
\end{align*}
where $(\lambda_k)$ is a sequence of strictly positive \emph{stepsizes}.
Rockafellar proved~\cite{Rock:ppa} that if \eqref{eq:g.mip} has a
solution and
\begin{align}
  \label{eq:et.ppm}
  \Norm{z_k-(\lambda_k T+I)^{-1}(z_{k-1})}\leq e_k,\;
   \sum_{k=1}^{\infty}e_k<\infty,\;\;\;
\inf \lambda_k>0,
\end{align}
then $(z_k)$ converges to a solution of \eqref{eq:g.mip}.

In each step of the PPM, computation of the proximal map
$(\lambda T+I)^{-1}z$ amounts to solving the \emph{proximal (sub)
  problem}
\begin{align*}
  0\in \lambda T(z_+)+z_+-z,
\end{align*}
a regularized inclusion problem which, although well posed,
is almost as hard as \eqref{eq:g.mip}.
From this fact stems the necessity of using approximations of the
proximal map, for example,  as prescribed in \eqref{eq:et.ppm}. 
Moreover, since each new iterate is, hopefully, just a better
approximation to the solution than the old one, if it were
compute with high accuracy, then the computational cost of each iteration would
be too high (or even prohibitive) and this would impair the overall
performance of the method (or even make it infeasible).

Unfortunately,
prescription~\eqref{eq:et.ppm} neither tells
which is the convenient error tolerance $e_k$ to be used in the
$k$-th iteration, nor  it guarantees that the next iterate will be
a better approximation than the current one. 

\section{Enlargements of maximal monotone operators and the Hybrid
Proximal Extragradient Method}

The Hybrid-Proximal Extragradient
Method~\cite{Sol-Sv:hy.ext,Sol-Sv:hy.proj} (hereafter HPE) is a
modification of the
PPM wherein \\ (a) the proximal subproblem, in each
iteration, is to be solved within a \emph{relative} error tolerance and
\\
(b) the update rule is modified so as to guarantee that the next iterate
is closer to the solution set by a quantifiable amount.
\\
An additional feature of (a) is that, in some sense, errors in the
inclusion on the proximal subproblems are allowed.
Recall that the \emph{$\varepsilon$-enlargement}~\cite{Bu-Iu-Sv:teps} of
a maximal monotone operator $T: H \tos H $ is
\begin{align}
  \label{eq:teps}
  T^{[\varepsilon]}(z)=\set{v\;|\; \inner{z-z'}{v-v'}\geq-\varepsilon
  \;\forall(z',v')\in T},\qquad 
  x\in H ,\;\varepsilon\geq0.
\end{align}

From now on, $T: H \tos H $ is a maximal monotone operator. The relaxed
HPE (r-HPE) method~\cite{Sol-Sv:hy.unif} for the monotone inclusion
problem~\eqref{eq:g.mip} proceed as follows: choose $z_0\in H $ and
$\sigma\in(0,1)$; for $k=1,2,\dots$
\begin{align}
  \label{eq:rhpe}
  \begin{aligned}
    &\text{compute }\tilde z_k,v_k,\varepsilon_k,{\lambda}_k>0\text{ such
      that} && v_k\in T^{[\varepsilon_k]}(\tilde z_k),\;\ \norm{{\lambda}_kv_k+\tilde
      z_k-z_{k-1}}^2+2{\lambda}_k\varepsilon_k
    \leq \sigma^2\norm{\tilde z_k-z_{k-1}}^2,\\
    &\text{choose }t_k\in(0,1]\text{ and
      set}&& z_k=z_{k-1}-t_k{\lambda}_kv_k.
  \end{aligned}
\end{align}
In practical applications, each problem has a particular structure which
may render feasible the computation of $\lambda_k$, $\tilde z_k$, $v_k$,
and $\varepsilon_k$ as above prescribed.
For example, $T$ may be Lipschitz continuous, it may be differentiable,
or it may be a sum of an operator which has a proximal map easily
computable with others with some of these properties.
Prescription for computing $\lambda_k$, $\tilde z_k$, $v_k$,
and $\varepsilon_k$ under each one of these assumptions were presented
in~\cite{Sol-Sv:hy.ext,
MR1682755,
MR1740961,
MR1953827,
MonSva10-3,
MonSva10-1,
MonSva11-1,
MonSva10-2%
}. 

An \emph{exact} PPM iteration for \eqref{eq:g.mip} is
$z_+=(\lambda T+I)^{-1}(z)$, where $z$ is the current iterate, $z_+$ is
the new iterate, $\lambda>0$ is the stepsize, and $I$ is the identity
map.
Computation of such a point $z_+$ is equivalent to solving, in the
variables $v,z_+$, the \emph{proximal inclusion-equation system}:
\[
 v\in T(z_+),\;\lambda v+z_+-z=0.
\]
Whence, the error criterion in~\eqref{eq:rhpe} relaxed both the
inclusion and the equality in the above inclusion-equation system.
The next lemma shows that an approximate solution of the proximal
inclusion-equation system satisfying that error criterion still
conveys useful information for solving \eqref{eq:g.mip}.

\begin{lemma}
  \label{lm:a}
  Take $z\in  H $; suppose that $\lambda>0$, $\sigma\in[0,1)$, $t\in[0,1]$,
  $\varepsilon>0$
  \begin{align}
    \label{eq:hpe.ec}
    v\in T^{[\varepsilon]}(\tilde z),\;\;
    \norm{\lambda v+\tilde z-z}^2+2\lambda \varepsilon \leq \sigma^2
    \norm{\tilde z-z}^2;
  \end{align}
  and define $z_+=z-t\lambda v$, $\gamma: H \to\R$,
  $\gamma(z')=\inner{z'-\tilde z}{v}-\varepsilon$. Then
  \begin{enumerate}
  \item
    \label{it:1-a}
    $(1-\sigma)\norm{\tilde z-z} \leq \norm{\lambda v} \leq
    (1+\sigma)\norm{\tilde z-z}$
    and $2\lambda\varepsilon\leq\sigma^2\norm{\tilde z-z}$;
  \item
    \label{it:2-a}
    $z_+=\argmin_{z'\in H } t\lambda\gamma(z')+\norm{z'-z}^2/2$;
  \item
    \label{it:3-a}
    $\min_{z'\in H }t\lambda\gamma(z')+\dfrac12\norm{z'-z}^2 \geq
   \dfrac12\Big(
  (1-\sigma^2)t\norm{\tilde z-z}^2+t(1-t)\norm{\lambda v}^2\Big)$;
 \item
   \label{it:4-a}
   for any $z^*\in T^{-1}(0)$,
   $\gamma(z^*)\leq 0$ and
   $\norm{z-z^*}^2 \geq \norm{z_+-z^*}^2+(1-\sigma^2)t\norm{\tilde z-z}^2
   +t(1-t)\norm{\lambda v}^2$;
 \item
   \label{it:5-a}
   for any $z^*\in T^{-1}(0)$,
   $\norm{z^*-\tilde z} \leq\norm{z^*-z}/\sqrt{1-\sigma^2}$ and
   $\norm{\tilde z-z} \leq \norm{z^*-z}/\sqrt{1-\sigma^2}$.
 \end{enumerate}
\end{lemma}

\begin{proof}[Proof of Lemma~\ref{lm:a}]
  Since $\lambda>0$ and $\varepsilon\geq 0$,
  $\norm{\lambda v+\tilde z-z} \leq \sigma\norm{\tilde z-z}$.
  Combining this inequality with triangle inequality
  we conclude that the two first inequalities in item~\ref{it:1-a} hold.
  The last inequality in item~\ref{it:1-a} follows trivially from
  the assumptions~\eqref{eq:hpe.ec}.
  Item~\ref{it:2-a} follows trivially from the definitions of $\gamma$
  and $z_+$.

  Direct use of item~\ref{it:2-a} and of the definitions of $z_+$ and
  $\gamma$ yields
  \begin{align*}
    \min_{z'\in{\R^p}}t\lambda\gamma(z')+\dfrac12\norm{z'-z}^2=
    \dfrac{1}{2}\Big(t\left[\norm{\tilde z-z}^2
    -\left(\norm{\lambda v+\tilde z-z}^2
    +2\lambda\varepsilon\right)\right]+t(1-t)\norm{\lambda v}^2\Big),
  \end{align*} which, combined with the inequality in \eqref{eq:hpe.ec} proves
  item~\ref{it:3-a}.

  The first inequality in item~\ref{it:4-a} follows from the inclusion
  $v\in T^{[\varepsilon]}(\tilde z)$, the definition of $\gamma$, and 
  the definition of $T^{[\varepsilon]}$  \eqref{eq:teps}, with $z'=z^*$
  and $v'=0$.
  Since $\lambda>0$, $t\geq0$, and $\gamma$ is affine, it follows from the
  first inequality in item~\ref{it:4-a}, item~\ref{it:2-a} and item~\ref{it:3-a}
  that
  \begin{align}
    \dfrac{1}{2}\norm{z^*-z}^2
    \geq    t\lambda\gamma(z^*)+\dfrac{1}{2}\norm{z^*-z}^2=
    \dfrac{1}{2}\norm{z^*-z_+}^2+t\gamma(z_+)+\dfrac{1}{2}\norm{z_+-z}^2
  \end{align}
  which, combined with items~\ref{it:2-a}
 and \ref{it:3-a} proves the second inequality in item
 \ref{it:4-a}.

 To prove the last item, define $\hat z=z-\lambda v$. Using item~\ref{it:4-a}
 with $t=1$, $z'=\hat z$ and the inequality \eqref{eq:hpe.ec} we conclude that
 \begin{align*}
   \norm{z^*-z}^2\geq\norm{z^*-\hat z}^2+(1-\sigma^2)\norm{\tilde z-z}^2,\quad
 \sigma\norm{\tilde z-z}\geq   \norm{\hat z-\tilde z} 
 \end{align*}
 Therefore, 
 \begin{align*}
   \norm{z^*-\tilde z} \leq\norm{z^*-\hat z}+\norm{\hat z-\tilde z}
   & \leq \norm{z^*-\hat z}+\sigma\norm{\tilde z-z}\\
   & \leq \sqrt{\norm{z^*-\hat z}^2+(1-\sigma^2)\norm{\tilde z-z}^2}\sqrt{1+\dfrac{\sigma^2}{1-\sigma^2}} \leq \dfrac{\norm{z^*-z}}{\sqrt{1-\sigma^2}}
 \end{align*}
 where the fist inequality follow from triangle inequality and the third
 from Cauchy-Schwarz inequality.
\end{proof}

In the next proposition we show that $z_k$ is closer than $z_{k-1}$ to
the solution set, with respect to the norm square, by a quantifiable
amount and derive some useful estimations.

\begin{proposition}
  \label{pr:a1}
  For any $k\geq 1$ and $x^*\in T^{-1}(0)$,
  \begin{enumerate}
  \item
    \label{it:1-b}
    $(1-\sigma)\norm{\tilde z_k-z_{k-1}} \leq \norm{\lambda_k v_k} \leq
    (1+\sigma)\norm{\tilde z_k-z_{k-1}}$
    and $2\lambda_k\varepsilon_k \leq \sigma^2 \norm{\tilde z_k-z_{k-1}}^2$;
  \item
    \label{it:2-b}
    $\norm{z^*-z_{k-1}}^2 \geq \norm{z^*-z_k}^2+t_k(1-\sigma^2)\norm{\tilde
      z_k-z_{k-1}}^2\geq\norm{z^*-z_{k-1}}^2$;
  \item
    \label{it:3-b}
    $\norm{z^*-z_0}^2 \geq \norm{z^*-z_k}^2+(1-\sigma^2)\sum_{j=1}^k
    t_j\norm{\tilde z_j-z_{j-1}}^2$;
  \item \label{it:4-b}
    $\norm{z^*-\tilde z_k}\leq\norm{z^*-z_{k-1}}/\sqrt{1-\sigma^2}$ 
    and $\norm{\tilde z_k-z_{k-1}} \leq \norm{z^*-z_{k-1}}/\sqrt{1-\sigma^2}$.
  \end{enumerate}
\end{proposition}

\begin{proof}
  Items \ref{it:1-b} and \ref{it:2-b} follow trivially from
  Lemma~\ref{lm:a}, items \ref{it:1-a} and \ref{it:4-a}, and the
  assumption $\sigma\in[0,1)$.
  Item~\ref{it:3-b} follows from item~\ref{it:2-b}.
  Item~\ref{it:4-b} follows from   Lemma~\ref{lm:a}, item \ref{it:5-a}.
\end{proof}

The aggregate stepsize $\Lambda_k$ and the ergodic sequences
$(\tilde{z}_k^a)$, $(\tilde v_k^a)$, and $(\varepsilon_k^a)$ associated
with the sequences $(\lambda_k)$, $(\tilde{z}_k)$,  $(v_k)$, and
$(\varepsilon_k)$ are, respectively,
\begin{align}
\label{eq:d.eg}
  \begin{aligned}
    &\Lambda_k:=\sum_{i=1} ^kt_i{\lambda}_i,\\
    &\tilde{z}_k^{\,a}:= \frac{1}{\;\Lambda_k}\;\sum_{i=1}^kt_i {\lambda}_i
    \tilde{z}_i, \quad v_k^{\,a}:= \frac{1}{\;\Lambda_k}\;\sum_{i=1}^k
    t_i{\lambda}_i v_i, \quad \varepsilon_k^{\,a}:=
    \frac{1}{\;\Lambda_k}\;\sum_{i=1}^k t_i{\lambda}_i (\varepsilon_i
    +\inner{\tilde{z}_i-\tilde{z}_k^{\,a}}{v_i-v_k^{\,a}}).
  \end{aligned}
\end{align}
The relevance of these ergodic sequences rests on the following theorem.

\begin{theorem}
  \label{th:a2}
  For any $k\geq 1$,
 $   v_k^a\in T^{\left[\varepsilon_k^a\right]}(\tilde z_k^a)$.
 Moreover, if $d_0$ is the distance from $z_0$ to
 $T^{-1}(0) \neq \emptyset$, then
 \begin{align*}
   \Norm{v_k^a}\leq \dfrac{2d_0}{\Lambda_k},\qquad
   \varepsilon_k^a \leq \dfrac{2d_0^2}{\Lambda_k\sqrt{1-\sigma^2}}
 \end{align*}
 for any $k \geq 1$
\end{theorem}

\begin{proof}[Proof of Theorem~\ref{th:a2}]
  The first part of the theorem follows from
  definitions~\eqref{eq:d.eg}, the inclusion in~\eqref{eq:rhpe},
  and the transportation formula for the
  $T^{[\varepsilon]}$~\cite[Theorem 3.11]{MR1716028}.

  To prove the second part of the Theorem, let $z^*$ be the projection
  of $z_0$ onto $T^{-1}(0)$. It follows from Proposition~\ref{pr:a1}
  item~\ref{it:2-a} that $\norm{z^*-z_k}\leq\norm{z^*-z_0}=d_0$ for any $k$.
  Theretofore,
  \begin{align}
    \label{eq:aux.a1}
    \norm{z_k-z_0}\leq 2d_0,\qquad \forall k\in\N.
  \end{align}
  Direct use of  
  the update rule for $z_k$ in \eqref{eq:rhpe} and of the 
  definition of $\Lambda_k$ and $v_k^a$ in \eqref{eq:d.eg} yields
  \begin{align}
    \label{eq:z0zk}
    z_0-z_k=\sum_{j=1}^kt_j\lambda_jv_k=\Lambda_kv_k^a.
  \end{align}
  The first inequality follows from the above equation and \eqref{eq:aux.a1}.

  Define, for $k=1,\dots$, the affine functions $\gamma_k,\Gamma_k:{\R^p}\to\R$,
  \begin{align}
    \label{eq:gammas}
    \gamma_k(z)=\inner{z-\tilde z_k}{v_k}-\varepsilon_k,
    \quad
    \Gamma_k(z)=\sum_{j=1}^k t_j\lambda_j\gamma_j(z).
  \end{align}
  We claim that for $k=1,2,\dots$
  \begin{align}
    \label{eq:new}
    z_k=\argmin_{z\in{\R^p}}\Gamma_k(z)+\dfrac12\norm{z-z_0}^2,
    \quad
    \min_{z\in{\R^p}}\Gamma_k(z)+\dfrac12\norm{z-z_0}^2
    \geq 0.     
  \end{align}
  The first above relation follow trivially from \eqref{eq:z0zk}. It
  follows from Lemma~\ref{lm:a}, items~\ref{it:2-a} and
  \ref{it:3-a} and the assumption $0\leq \sigma\leq 1$ that          the second
  relation in \eqref{eq:new} holds for $k=1$.
  If the inequality in \eqref{eq:new} holds for $k$, as
  $\Gamma_{k+1}=\Gamma_k+t_{k+1}\lambda_{k+1}\gamma_{k+1}$, it follows
  again from Lemma~\ref{lm:a} items~\ref{it:2-a} and \ref{it:3-a}, the
  assumption $0\leq \sigma\leq 1$, and the first relation in
  \eqref{eq:new} that this inequality holds for $k+1$.

  It follows from~\eqref{eq:new} and definitions~\eqref{eq:gammas} that 
  \begin{align*}
    \Gamma_k(\tilde z_k^a)+\dfrac{1}{2}\norm{\tilde z_k^a-z_0}^2
    &=\dfrac{1}{2}\norm{\tilde z_k^a-z_k}^2
      +\Gamma_k(z_k)+\dfrac{1}{2}\norm{z_k-z_0}^2
      \geq\dfrac12\norm{\tilde z_k^a-z_k}^2
  \end{align*}
  Direct use of the transportation formula for the
  $T^{[\varepsilon]}$~\cite[Theorem 3.11]{MR1716028},
  \eqref{eq:d.eg}, and \eqref{eq:gammas} shows that
  $-\Gamma_k(\tilde z_k^a)=\Lambda_k\varepsilon_k^a$. Therefore
  \begin{align*}
    2\Lambda_k\varepsilon_k
    \leq \norm{\tilde z_k^a-z_0}^2-\norm{\tilde z_k^a-z_k}^2
    &=2\inner{\tilde z_k^a-z_0}{z_k-z_0}-\norm{z_k-z_0}^2\\
    & \leq 2d_0\left(1+\dfrac{1}{\sqrt{1-\sigma^2}}\right)
      \norm{z_k-z_0}-\norm{z_k-z_0}^2  \leq\dfrac{4d_0^2}{\sqrt{1-\sigma^2}}
  \end{align*}
\end{proof}

Next we analyze the pointwise and ergodic complexities
 of the r-HPE when 
the large-step condition, introduced
in~\cite{MonSva10-3,MonSva11-1}, is satisfied and the relaxation parameters
$t_k$ are bounded away from zero.

\begin{theorem}
  \label{lm:rhpe2}
  Let $d_0$ be the distance from $z_0$ to $T^{-1}(0)\neq\emptyset$.
  If %
  for any $k\geq1$,
\begin{align}
  \label{eq:lsc}
  \lambda_k\norm{\tilde z_k-z_{k-1}} \geq
  \eta>0,\qquad t_k\geq \tau>0
\end{align}
 then, for any $k\geq1$,
  \begin{enumerate}
  \item\label{it:rhpe2-1}
    there exists $i$, $1 \leq i \leq k$, such that
    \begin{align*}
      \norm{v_i}\leq\dfrac{d_0^2}{\eta (1-\sigma) k \tau},\qquad
      \varepsilon_i\leq
      \dfrac{\sigma^2}{2\eta}\dfrac{d_0^3}{\left((1-\sigma^2)k\tau\right)^{3/2}};
    \end{align*}
  \item\label{it:rhpe2-2}
    $v_k^a\in T^{[\varepsilon_k^a]}(\tilde z_k^a)$,
    \begin{align*}
      \norm{v_k^a}\leq \dfrac{2d_0^2}{(\tau k)^{3/2}\eta\sqrt{1-\sigma^2}},
      \qquad
\varepsilon_k^a\leq
 \dfrac{2d_0^3}{(\tau k)^{3/2}\eta(1-\sigma^2)}.
    \end{align*}
  \end{enumerate}
\end{theorem}
\begin{proof}
  It follows from Proposition~\ref{pr:a1},
  item~\ref{it:3-b}, that there exists $1 \leq i \leq k$ such that
  \begin{align*}
    \norm{\tilde z_i-z_{i-1}}\leq \dfrac{d_0}{\sqrt{(1-\sigma^2)\tau k}}
  \end{align*}
  It follows from  the first part of Proposition~\ref{pr:a1} and
  \eqref{eq:lsc} that, in particular for such an $i$,
  \begin{align*}
    \norm{v_i}\leq \dfrac{(1+\sigma)\norm{\tilde z_i-z_{i-1}}}{\lambda_i},
    \quad
    \varepsilon_i \leq\dfrac{\sigma^2\norm{\tilde z_i-z_{i-1}}^2}{2\lambda_i},
    \quad
    \dfrac{1}{\lambda_i} \leq \dfrac{\norm{\tilde z_i-z_{i-1}}}{\eta}.
  \end{align*}
  Item~\ref{it:rhpe2-1} follows from the above inequalities.

  It follows from \eqref{eq:lsc} and Proposition~\ref{pr:a1}
  item~\ref{it:3-b} that
  \begin{align*}
    \sum_{j=1}^k \tau \dfrac{\eta^2}{\lambda_j^2}
    \leq \sum_{j=1}^k t_j\norm{\tilde z_j-z_{j-1}}^2
    \leq \dfrac{d_0^2}{1-\sigma^2}.
  \end{align*}
  Using this result and Lemma~\ref{lm:ti} we conclude that 
  \begin{align*}
    \sum_{j=1}^k \lambda_j \geq k^{3/2}\left(
    \dfrac{d_0^2}{\tau\eta^2(1-\sigma^2)}\right)^{-1/2}
    \quad \text{and}
    \quad \Lambda_k\geq (\tau k)^{3/2}\dfrac{\eta\sqrt{1-\sigma^2}}{d_0},
  \end{align*}
  where the second inequality follows from \eqref{eq:d.eg} and the assumption
  $t_j\geq \tau$ for all $j$.
  Item~\ref{it:rhpe2-2} follows from the second above inequality and
  Theorem~\ref{th:a2}.
\end{proof}

\appendix

\section{Auxiliary results}
\label{ap:ar}

\begin{lemma}
  \label{lm:ti}
  If $\alpha_i>0$ for $i=1,\ldots,m$ and $\sum_{i=1}^m\alpha_i^{-2}\leq C$ then
  \(
   \sum_{i=1}^m\alpha_i\geq m^{3/2}/C^{1/2}.
\)
\end{lemma}

\begin{proof}
Take $\alpha\in\R^m_{++}$  such that
  $\sum_{i=1}^m \alpha_i^{-2}\leq C$ and let $\bar\alpha=\sum_{i=1}^m\alpha_i/m$.
  As $t^{-2}$ is convex for $t>0$,
  \begin{align*}
    \dfrac{1}{{\bar\alpha}^2}\leq\dfrac{1}{m}\sum_{i=1}^k\dfrac{1}{\alpha_i^2}
    \leq \dfrac{C}{m};
  \end{align*}
  therefore, $\sqrt{m/C}\leq \bar \alpha$. To end the proof, use the definition of $\bar\alpha$.
\end{proof}

The next result was proved in \cite[Corollary 1]{error-bound}
\begin{lemma}
  \label{lm:d}
  If $T: H \tos H $ is maximal monotone, $z\in H $ and
  $\tilde v\in T^{[\varepsilon]}(\tilde z)$, then
  \begin{align*}
    \norm{\lambda\tilde v+\tilde z-z}^2+2\lambda\varepsilon
    \geq \Norm{\tilde z-(\lambda T+I)^{-1}z}^2
    +\Norm{\lambda^{-1}\bigg((\lambda T+I)^{-1}z-z\bigg)}^2.
  \end{align*}
\end{lemma}


\def\cprime{$'$}

\end{document}